\documentclass[12pt]{amsart}
\usepackage{amsthm,amsmath,amsxtra,amscd,amssymb,xypic,xparse,xspace,mathtools}
\usepackage[all]{xy}
\usepackage{tikz-cd}
\numberwithin{equation}{section}
\usepackage{bm}
\usepackage{hyperref}
\usepackage{cleveref}
\usepackage{tabularx}
\usepackage{booktabs}
\usepackage{caption,subcaption}
\usepackage{enumitem}  
\usepackage{mathrsfs}

\theoremstyle{plain}
\newtheorem{thm}{Theorem}[section]

\newtheorem{prop}[thm]{Proposition}
\newtheorem{cor}[thm]{Corollary}

\newtheorem{remark}[thm]{Remark}

\theoremstyle{definition}
\newtheorem{definition}[thm]{Definition}

\newtheoremstyle{colon}%
{}
{}
{\itshape}
{}
{\bfseries}
{:}
{ }
{}

\theoremstyle{colon}

\def\={\;=\;}  \def\+{\,+\,}

\newcommand{\GL}{\operatorname{GL}}

\newcommand{\rank}{\operatorname{rank}}

\newcommand{\calM}{{\mathcal M}}
\newcommand{\calO}{{\mathcal O}}

\newcommand{\calQ}{{\mathcal Q}}

\newcommand{\CC}{{\mathbb{C}}}

\newcommand{\PP}{{\mathbb{P}}}

\newcommand{\RR}{{\mathbb{R}}}
\newcommand{\ZZ}{{\mathbb{Z}}}

\newcommand{\poles}{\underline{p}}
\newcommand{\zeroes}{\underline{z}}

\newcommand{\gon}{\operatorname{gon}}

\newcommand{\lG}{{\Gamma}}     

\newcommand{\hor}{{\mathrm {hor}}}
\newcommand{\ver}{{\mathrm {ver}}}

\DeclareDocumentCommand{\Ehor}{O{\lG}}{E^{\hor}(#1)}   
\DeclareDocumentCommand{\Ever}{O{\lG}}{E^{\ver}(#1)}  
\DeclareDocumentCommand{\Ehori}{O{i} O{\lG} }{E^{\hor}_{(#1)}(#2)}   

\DeclareDocumentCommand{\abEhor}{O{\lG}}{\widetilde{E^{\hor}}(#1)}

\DeclareDocumentCommand{\relhomZ}{O{X} O{\zeroes} O{\poles}}{H_1(#1\setminus #3,#2;\ZZ)}
\DeclareDocumentCommand{\relhomC}{O{X} O{\zeroes} O{\poles}}{H_1(#1\setminus #3,#2;\CC)}

\usepackage{color}

\newcommand\TM{Teichm\"uller\xspace}

\title[Totally geodesic subvarieties and linear systems]{ Totally geodesic subvarieties of  the moduli space of curves and  linear systems}
\author{Frederik Benirschke}

\begin{document}

\maketitle

\begin{abstract}

We construct a linear system on a general curve in a totally geodesic subvariety of the moduli space of curves. As a consequence, we obtain rank bounds for totally geodesic subvarieties of dimension at least two. Furthermore, we classify totally geodesic subvarieties of dimension at least two in strata with at most two zeros.

\end{abstract}

\section{Introduction}

Let $\calM_{g,n}$ be the moduli space of genus $g$ curves with $n$ marked points.
 A {\em totally geodesic subvariety} (for the \TM metric)  of $\calM_{g,n}$ is an algebraic subvariety $M\subseteq \calM_{g,n}$ such that any \TM geodesic passing through a general point and tangent to $M$ is contained in $M$. Totally geodesic subvarieties are closely related to $\GL(2,\RR)$-orbit closures in strata of quadratic differentials. Let $QM$ be the subbundle of the bundle of quadratic differentials $\calQ_{\calM_{g,n}}$ containing pairs $(X,q)$ where $X\in M$ and $q$ is a quadratic differential generating a \TM geodesic contained in $M$.
The bundle $QM$ is stratified by the order of zeros and poles, and we let $QM(\mu)\subseteq \calQ(\mu)$ be the stratum of maximal dimension.
Then $QM(\mu)$ is a $\GL(2,\RR)$-orbit closure of dimension\[
\dim QM = \dim QM(\mu) = 2\dim M
\]
in some stratum of quadratic differentials $\calQ(\mu)\subseteq \calQ_{\calM_{g,n}}$.

In this paper, we mostly work with the orbit closure $QM(\mu)$ instead of the totally geodesic subvariety $M$. The following definition allows us to switch between the two points of view. 
\begin{definition}
Let $N\subseteq\calQ(\mu)$ be an orbit closure in a stratum of quadratic differentials.
We say $N$ is a {\em totally geodesic orbit closure} if 
\[
\dim N = 2\dim \pi(N),
\]
where $\pi:\calQ(\mu)\to\calM_{g,n}$ is the forgetful map.
\end{definition}
There is a $1$-to-$1$ correspondence between totally geodesic subvarieties and totally geodesic orbit closures given by 
\[
N = QM(\mu),\, M = \overline{\pi(N)};
\] see \cite{goujard,wright-totally-geodesic}.
The rank of a totally geodesic orbit closure is  \[\rank(N):= \tfrac{\dim N}{2}.\footnote{By taking square roots of quadratic differentials in $N$, one obtains an orbit closure $N'$ in a stratum of Abelian differentials. The rank of $N'$ is by definition half of the dimension of the projection of the tangent space $p(T_{(X,\omega)}N')\subseteq H^1(X,Z(\omega),\CC)$ at a generic point $(X,\omega)\in N'$.
By \cite[Thm 1.3]{wright-totally-geodesic} the rank of $N'$ is $\tfrac{\dim(N)}{2}$. Thus with our definition $\rank(N)=\rank(N').$
}
\]

Our first result is a rank bound for totally geodesic orbit closures, depending only on the number of simple zeros in the partition $\mu$.
Any partition $\mu $ can be written as \[
\mu = (-1^n,1^m, b_1,\ldots,b_k),\, b_i>1, i=1,\ldots k.
\]We call $m$ the {\em number of simple zeros}.
\begin{thm}\label{thm:bound} Let $\calQ(\mu)$ be a stratum of quadratic differentials in genus $g$ with $m$ simple zeros.
Suppose $N\subseteq \calQ(\mu)$ is a totally geodesic orbit closure of $\rank(N)\geq 2.$
Then 

\[
\rank(N) \leq \begin{cases} m-g+1 & \text{ if $m\geq 2g-1$},\\
\tfrac{m}{2}+1  &\text{ if $m\leq 2g-1$}.
\end{cases}
\]
In particular, if $\rank(N)= m+1$, then  $g(X) = 0$ and $\rank(N)\leq n-3$.

\end{thm}

The special case $m=0$ rules out the existence of totally geodesic subvarieties of rank at least $2$ in strata with only higher-order zeros.

\begin{cor}
There do not exist totally geodesic orbit closures of rank at least $2$ in strata  $\calQ(-1^n,b_1,\ldots,b_k)$ with $b_i>1$ for $ i=1,\ldots,k$.
\end{cor}

We prove \Cref{thm:bound} by constructing a linear system of degree $m$ and dimension $\rank(N)$ and applying results from the theory of special divisors; see \Cref{prop:exist} for a more precise statement. In the special case of a totally geodesic surface $N$, we obtain an upper bound on the gonality of curves in $\pi(N)$; see \Cref{prop:gon}.

 A basic construction for producing totally geodesic subvarieties is covering constructions,  which are obtained by pulling back differentials along branched coverings. The resulting orbit closures are called {\em loci of covers}.
 {\em Primitive} totally geodesic subvarieties are the totally geodesic subvarieties that do not arise as covering construction. It is an important open problem to classify all primitive totally geodesic subvarieties.
 Mirzakhani initially conjectured that all totally geodesic subvarieties of dimension at least $2$ are covering constructions. Still, three examples of primitive totally geodesic surfaces were recently found by McMullen, Mukamel, Wright, and Eskin \cite{mmw, emmw}.

We rule out the existence of primitive, totally geodesic orbit closures in strata of quadratic differentials with at most two zeros. I.e., strata of the form
$\calQ(-1^n,a,b), a,b \geq 0$ or $\calQ(-1^n,a), a\geq 0$.

\begin{thm} \label{thm:main}
Let $\calQ(\mu)$ be a stratum of quadratic differentials with at most two zeros.
Suppose $N\subseteq \calQ(\mu)$ is a totally geodesic orbit closure of rank at least $2$.
Then $N$ is one of the following strata 
\[
\calQ(-1^5,1), \calQ(-1^6,1^2), \calQ(-1^2,1^2).
\]

 In particular, there do not exist primitive totally geodesic orbit closures in $\calQ(\mu)$.
\end{thm}

The known examples of primitive totally geodesic orbit closures from  \cite{emmw,mmw}  are totally geodesic surfaces in the strata
\[
\calQ(-1^3,1^3), \calQ(-1^4,1^4), \calQ(-1,1^5).
\]

 Thus,  \Cref{thm:main} cannot be extended to strata of quadratic differentials with more than two zeros without further assumptions.

\subsection*{Acknowledgements}
We thank Paul Apisa, Alex Wright, Carlos Serv{\`a}n, Dawei Chen and Aaron Calderon for helpful conversations.

\section{Linear systems of totally geodesic orbit closures}

Let $M$ be a totally geodesic subvariety of $\calM_{g,n}$. Our main tool is the construction of a linear series for a general curve in $M$. Recall that a linear series $\mathscr{D}$ on a curve $X$ is a linear subspace $V\subseteq H^0(\calO_X(D))$ for a divisor $D$ on $X$. A linear series of {\em dimension} $ r+1 := \dim(V) $ and {\em degree} $d := \deg(D)$ is called a $\mathfrak{g}^r_d.$ We refer the reader to \cite[Chap. I + III]{acgh} for an introduction to the theory of linear series. 
The construction of the linear series starts with the following observation,  which is a variation of  \cite[Prop. 3.1]{FredCarlos}.
Recall that $QM\subseteq \calQ_{\calM_{g,n}}$ is the bundle of all quadratic differentials generating \TM geodesics that are contained in $M$.
 \begin{prop}
Suppose $M\subseteq\calM_{g,n}$ is a totally geodesic subvariety of dimension $d$  and $X\in M^{reg}$ a regular point. Then, the fiber of $QM$ over $X$ is a linear subspace of dimension $d$.
 \end{prop}

 \begin{proof}
 Consider the map \[
\phi: Q_XM\hookrightarrow \calQ_{\calM_{g,n},X}  = T^*_X{\calM_{g,n}}\to T^*_XM,
 \]
 where $Q_XM$ is the fiber of $QM$ over $X$ and $T^*_XM$ is the cotangent space to $M$.

 Following the argument in  \cite[Prop. 3.1]{FredCarlos}, the map $\phi$ is injective, proper, homogenous.  The homogeneity follows from $\GL(2,\RR)$-invariance of $QM$, while injectivity and properness follow from $\phi$ being norm-preserving for the restriction of the $L^1$-norm on $Q_XM$ and the quotient $L^1$-norm on $T^*_XM$. Since $QM$ is an algebraic variety, $\phi$ is also holomorphic.

 Let \[
 U:=(Q_XM)^{reg},\, Z:= Q_XM \setminus U.
 \]
 Since $\phi$ is proper, $\phi(Z)$ is a proper closed subvariety and $T^*_XM\setminus \phi(Z)$ is connected.  By invariance of domain, $\phi_{|U}$ is a homeomorphism onto  $T^*_XM\setminus \phi(Z)$ and hence $\phi$ is surjective.
 We conclude that $\phi$ is a homeomorphism. Restricted to $T^*_XM\setminus \phi(Z)$ the inverse $\phi^{-1}$ is holomorphic and by the Riemann extension theorem, $\phi^{-1}$ is holomorphic.
 Additionally, $\phi^{-1}$ is homogeneous. Thus, $\phi^{-1}$ is homogeneous and differentiable at the origin and hence linear.
 It follows that \[
 Q_XM = \phi^{-1}(T^*_XM)
 \]
 is a linear subspace of dimension \[
 \dim T^*_XM = \dim M=d.
 \]
 \end{proof}

We now switch from totally geodesic subvarieties to totally geodesic orbit closures.
Let $N\subseteq \calQ(\mu)$ be a totally geodesic orbit closure of rank at least $2$ in a stratum of quadratic differentials.
 Suppose  \[
 (X,q)\in N\subseteq\calQ(\mu).
 \]For the rest of this section, we write 
 \[
 \begin{gathered}
\mu = (-1^n,1^m, b_1,\ldots,b_k), b_i>1,\\
g= g(X),\, r = \rank(N)-1,\\
 (q) = B + \sum_{i=1}^m x_i - P,\\
  B= \sum_{i=1}^k b_i y_i, P = \sum_{i=1}^n p_i \text{ for some points $x_i, p_i,y_i\in X$}.
 \end{gathered}
 \]

Let \[
  \mathscr{D}_X := \{ P + (q),\,\, (X,q)\in N\}
  \]
 be the associated  $\mathfrak{g}^r_{4g-4+n}$. The next result estimates the base locus of  $\mathscr{D}_X$.

\begin{prop}\label{prop:exist}
Suppose $N\subseteq\calQ(\mu)$ is a totally geodesic orbit closure with  $\rank(N)\geq 2$ and $X\in \pi(N)$ is generic. 
Then $B$ is contained in the base locus of $\mathscr{D}_X$ and thus
\[
\mathscr{D}_X = B + \mathscr{D}'_X,
\]
where $\mathscr{D}'_X$ is a $\mathfrak{g}^r_m$.
In particular, a generic curve $X\in \pi(N)\subseteq\calM_{g,n}$ admits a $\mathfrak{g}^r_m$.
\end{prop}

\begin{proof}
Let $B'$ be the base locus of $\mathscr{D}_X$. The goal is to show that $B'\geq B$.
A general divisor in $\mathscr{D}_X$ is of the form
\[
P + (q) = B' + \sum_{i=1}^{4g-4+n-\deg(B')} z_i
\]
where $z_i\in X$ are pairwise distinct and not contained in $B'$.
On the other hand, a generic quadratic differential $q\in V$ is contained in 
\[
\calQ(-1^n,1^m, b_1,\ldots,b_k), b_i>1.
\]
In particular each point $z_i$ is one of the $m$ simple zeros. Hence $B\leq B'$.

\end{proof}

 \Cref{thm:bound} now follows from the existence of a $\mathfrak{g}^r_m$, together with Riemann-Roch and Clifford's theorem.

\begin{proof}[Proof of \Cref{thm:bound}]
By \Cref{prop:exist}, there exists a $\mathfrak{g}^r_m$ on a generic curve $X$ in $N$. It follows from Riemann-Roch that if $m\geq 2g-1$, one has $r\leq m-g$.
The case $r\leq 2g-1$ follows from Clifford's theorem \cite[III.1]{acgh}.
 In case of $r=m$ there exists a divisor $D$ on $X$ with $h^0(D) = \deg(D)+1$. Hence, there exists a degree one map from $X$ to $\PP^1$. Thus $g(X) =0$.
 \end{proof}

\begin{remark}
Another consequence of Clifford's theorem is that 
if 
\[
m<2g-2,\, 2r = m,
\] then $\pi(N)\subseteq \calM_{g,n}$ is contained in the hyperelliptic locus.
We will see in \Cref{rem:hyp} that the same conclusion holds in the case  $m <g, \, 2r = m$.
\end{remark}

We now turn to the  proof of  \Cref{thm:main}.

\begin{proof}[Proof of \Cref{thm:main}]
First, assume that $\calQ(-1^n,a)$ is a stratum with only one zero. By \Cref{thm:bound},  $m=1$ and thus $N$ has to be one of the following strata $\calQ(-1^5,1), \calQ(-1,1)$. The former has rank $2$, while the second stratum is empty.

Next, we adress strata   $\calQ(-1^n,a, b)$ with two zeros. The first step is to show that $N$ is a locus of covers in a stratum in genus zero or one.  

By \Cref{thm:bound}, the case $m=0$ is impossible, and $m=1$ is only possible if $g=0$. If $m=1,g=0$, it follows from \cite[Thm. 1.5]{Apisa} that $N$ is a locus of covers.

The remaining case is $m=2$. This can only occur in the two strata 
\[
\calQ(-1^6,1^2), \calQ(-1^2,1^2).
\]
The first case is a stratum in genus zero with two zeros, in which case we can apply  \cite[Thm. 1.5]{Apisa} again to conclude that $N$ is a locus of covers. The second case would lead to a totally geodesic orbit closure $N$ with $\pi(N)\subseteq \calM_{1,2}$.
For dimension reasons this is only possible if $\pi(N) = \calM_{1,2}$ and $N$ is the generic stratum in genus $1.$

So far, we have shown that either $g=0$ and $N$ is a locus of covers or $g=1$ and $N=\calQ(-1^2,1^2)$. It remains to analyze the cases with $g=0$.
By \Cref{thm:bound} it follows $\rank(N)\leq 3$. If $\rank(N)=3$ this implies $m=2, \dim N=6$ and thus $N= \calQ(-1^6,1^2)$.
The remaining case is $\rank(N)=2$ and $\calQ(\mu)= (-1^{a+5},1,a)$.
Since $N$ is a full locus of covers, it has to be a cover of
a $4$-dimensional stratum in genus zero. The only possibilities are 
\[
\calQ(-1^4,0^2), \calQ(-1^5,1).
\]
The first one is not possible since, in this case, the orbit closure has non-zero rel, and totally geodesic orbit closures always have zero rel, see \cite[Thm. 1.3]{wright-totally-geodesic}.
In the second case, we will reach a contradiction from Riemann-Hurwitz.
The zero of order $a$ has to lie over the simple zero in $\calQ(-1^5,1)$.
In particular a pullback under a degree $d$ map has at least $5d-(2+k)$ simple poles, where $k$ is the number of simple branch points lying over the simple pole.   The additional $2$ corresponds to a potential ramification point of multiplicity $3$, in case the simple zero in $\calQ(-1^{a+5},1,a)$ lies over a simple pole. By Riemann-Hurwitz we have $k\leq 2d-2$.
Thus we obtain
\[
 5 \geq  5d-2-k \geq 3d,
\]
which is impossible for $d\geq 2.$

\end{proof}

It follows from \Cref{thm:bound} that $2r\leq m$ if $m\leq 2g-1$.
The following result gives further restrictions in the range $2r\leq m \leq 3r-2$.

\begin{cor}
Let $\calQ(\mu)$ be a stratum of quadratic differentials with $m$ zeros and $N\subseteq \calQ(\mu)$ a totally geodesic orbit closure of rank $r+1\geq 2$.
Let $\alpha = m-2r$ and assume $0\leq \alpha \leq r-2$.
    Then, one of the following is true
    \begin{itemize}
        \item  $g\leq r+2\alpha+1$ or,
        \item a generic curve $X$ in $\pi(N)$ is a double cover of a curve of genus at most $\tfrac{\alpha}{2}$.
    \end{itemize}
\end{cor}

\begin{proof}
By \cite[Chap. 3, Exc. B-7]{acgh} the existence of a \[
\mathfrak{g}^{r}_{2r+\alpha},\, 0\leq \alpha\leq r-2
\]
implies that either  $g\leq r+2\alpha+1$ or that a generic curve $X$ in $N$ is a double cover of a curve of genus $g'\leq \tfrac{\alpha}{2}$.
\end{proof}

\begin{remark}\label{rem:hyp}
For example, the statement for $\alpha =0$ in the above proposition says that  if $N$ is a totally geodesic orbit closure with \[
2\leq \rank(N)= \tfrac{m}{2}+1 < g,
\]then $\pi(N)$ is in the hyperelliptic locus.
\end{remark}

\begin{remark}
Even if $m$ is not in the range between $2r$ and $3r-2$, the linear system $\mathscr{D}_X$ induces  a map $\phi:X\to \PP^d, d\leq m$, after removing the base locus. There are two cases. Either $\phi$ is birational, in which case one can bound the genus of $X$ only in terms of $r$ and $m$ using Castelnuovo's bound (see \cite[Chapt. III.2]{acgh}). The second case is that $X$ covers another curve with a degree of at least two. It would be interesting to find a numerical criterion for when the totally geodesic orbit closure is a locus of covers in the second case.
\end{remark}

\subsection*{Gonality bounds for totally geodesic surfaces}
Recall that the gonality $\gon(X)$ of a curve $X$ is the smallest degree of a non-constant holomorphic map to $\PP^1$.
The gonality of a general curve $X$ of genus $g$ is 
\[
\gon(X) = \lfloor \tfrac{g+3}{2}\rfloor.
\]
In \cite{bud}, Bud showed that in many strata of quadratic differentials, for example, if the partition has only positive entries, a generic curve still has gonality $ \lfloor \tfrac{g+3}{2}\rfloor$.
On the other hand, we obtain gonality bounds for totally geodesic surfaces only in terms of the number of simple poles of the partition $\mu$. 

\begin{cor}\label{prop:gon}
Suppose $N\subseteq \calQ(\mu)$ is a totally geodesic orbit closure of rank $2$ and $(X,\omega) \in N$. Then 
\[
\gon(X) \leq m,
\]
where $m$ is the number of simple zeros of $\mu$.
In particular, if $m=1$, then $g(X)=0$ and, if $m=2$, then $X$ is hyperelliptic.
\end{cor}
\begin{proof}
If $\rank(N)=2$,  the linear system $\mathscr{D}$ from \Cref{prop:exist} is a $\mathfrak{g}^1_m$ and hence \[
\gon(X)\leq m.
\]
\end{proof}

\begin{remark}
For any rank $2$ orbit closure $M$ in a stratum of Abelian differentials, one can construct a map to $\PP^1$ similarly. Let $(X,\omega)\!~\in~\!\pi(M)$. Consider the projection of the tangent space of $M$ to absolute cohomology $H^1(X;\CC)$. The intersection with $H^{1,0}(X)$ is a $2$-dimensional subspace and hence defines a map to $\PP^1$ after removing the base locus. If one can prove a lower bound on the size of the base locus, similar to the bound we have for totally geodesic orbit closures in terms of the number of simple zeros, then one would also obtain gonality bounds for curves in $M$.
\end{remark}

\end{document}